\documentclass[a4paper,11pt,reqno]{amsart}

\usepackage{tikz}
\usepackage{color}
\definecolor{Red}{rgb}{0.3,0.3,0.9}
\usepackage{multicol}

\usepackage{graphicx,colordvi}

 \newcounter{zum}
 
\usepackage{amssymb, amsmath}

\usepackage[active]{srcltx}

\begin{document}




\begin{center}
{\bf A special modulus of continuity and the $K-$functional}\\
{Nadezhda Dolmatova}\\
{Institute of Mathematics, University of Wroc{\l}aw, pl. Grunwaldzki 2/4, 
50-384 Wroc{\l}aw, Poland }
\end{center}

\section*{\bf Abstract}
 We consider the questions connected with the approximation of a real continuous 1 - periodic functions and give a new proof of the equivalence of the special Boman--Shapiro modulus of continuity with Peetre's $K-$functional. We also prove  Jackson's inequality for the approximation by trigonometric polynomials.


keywords: Modulus of continuity, $K-$functional, Jackson's theorem 

[2010]  42A10, 41A17, 41A50 



\section{\bf Introduction}
\setcounter{equation}{0}

Denote by $ C (\mathbb{T}), \mathbb{T} = \mathbb{R} / \mathbb{Z} $ the space of real continuous $1-$periodic functions $ f $ with the uniform norm
\begin{equation*}
\| f \| = \sup_{u \in \mathbb{T}} |f (u) |.
\end{equation*}
The derivative operator is denoted by the symbol $ D $, and the space of functions $f$ 
with $ D^2f \in C(\mathbb{T}) $ will be denoted by 
$C^2 (\mathbb{T}). $

Let $ L(\mathbb{T})$ be the space of measurable, integrable functions with norm
\begin{equation*}
\| f \|_1 = \int_{\mathbb{T}} |f(u)| \,du
\end{equation*}
and let $T_{n-1} $ be the set of real trigonometric $ 1-$ periodic
  polynomials $\tau $ of degree at most $n-1$:
\begin{equation*}
 \tau(t):=\sum_{j=-n+1}^{n-1} \alpha_j \exp (2\pi i j t), \qquad \alpha_j = \overline \alpha_{-j}.
\end{equation*}
 
For $ f \in C(\mathbb{T)} $, we denote by $ E_{n-1}(f) $ the value of the best approximation of $f$ by real trigonometric polynomials of degree at most $n-1$
\begin{equation*}
E_{n-1}(f):=\inf_{\tau \in T_{n-1}} \| f - \tau \|.
\end{equation*}

We will use the convolution of periodic functions $ f $ with positive functions $ g, $ with
finite support. In this case, the convolution can be understood in the following sense:

\begin{equation*}
(f*g)(t):= \int_{\mathbb R} f(u) \, g(t-u) \, du.
\end{equation*}

We denote by $ \chi_h^k,\ k = 1,2, \dots $ the convolution powers of the normalized characteristic function of the interval $ (-h/2,h/2),\ h>0 $:
\begin{equation*}
\chi_h^k:=\chi_h^{k-1}*\chi_h, \quad \chi_h(t):=
\begin{cases}\label{a}
\frac{1}{h}&t\in(-h/2,h/2),\\
0,&t\not\in(-h/2,h/2).
\end{cases}
\end{equation*}
In particular,
\begin{equation*}
\chi_h^2(t)=\begin{cases}
\frac{1}{h}\left(1-\frac{|t|}{h}\right),&t\in (-h,h),\\
0,& t\not\in (-h,h).
\end{cases}
\end{equation*}
The functions $\chi_h^k$ are  the cardinal $B-$splines with support $[-kh/2, kh/2]$ and 
$\| \chi_h^k \|_1 = 1$.
 
We will use the following moduli of continuity (see \cite{fks09,bdk12,bks13})
\begin{align*}\label{wh}
W_2(f,\chi_h^k)&:=\| f - f*\chi_h^k \|,\\W_2^*(f,\chi_h^k)&:=\sup_{0<u\leq h} W_2(f,\chi_u^k).
\end{align*}
\noindent
They are special cases of the Boman--Shapiro  moduli of continuity (see \cite{sh68,bsh71,b80}).

This paper is the continuation of \cite{bdk12}. The main result of \cite{bdk12}  is the following Jackson inequality for the uniform approximation of continuous $1-$periodic functions by trigonometric polynomials:
\vskip0.2cm
\textit{
Let $f$ be a continuous $1-$periodic function and $n\in \mathbb{N}$, $h=\alpha/(2n)$, $\alpha > 2/\pi$. Then the following inequality holds  
\begin{equation}\label{aa}
E_{n-1}(f)\leq  (\sec{1/\alpha} + \tan{1/\alpha})\, W_2(f,\chi_h).
\end{equation}
The estimate is exact for $ \alpha = 1,3, \dots $.
}
\vskip0.2cm
In  \cite{bdk12}  the following sharp Bernstein--Nikolsky--Stechkin inequality for $\tau \in T_n$ was also obtained:
\vskip0.2cm
\textit{
Let $\tau$ be a real trigonometric $ 1-$periodic
  polynomial of degree at most $n-1$ for $n\in \mathbb{N}$, and suppose $h\in(0,1/n]$. Then
\begin{equation}\label{b}
\| D^2 \tau \| \le  (2\pi n)^2 \ W_2(c_n, \chi_h)^{-1} \, W_2(\tau, \chi_h), \quad c_n(t):=\cos(2\pi n t).
\end{equation} 
}

The Jackson inequality \eqref{aa} and the Bernstein--Nikolsky--Stechkin estimate \eqref{b} allowed us to prove the equivalence of a special modulus of continuity and the second Peetre's $ K-$functional \cite{bdk12}:
\vskip0.1cm
\textit{
Let $h\in(0,1]$. Then
\begin{equation}\label{k}
1/4 K_2(f, h/(4\sqrt{6})) \le W_2 (f, \chi_h)\le 4 K_2(f, h/(4\sqrt{6})).
\end{equation}}

The  equivalence of moduli of this type and the $K-$ functional is known (see \cite{di93} and \cite{t05}). Here we give a new form of this equivalence with the calculation of the constants.
We present a new simple proof of the estimates of the type  \eqref{k} with
better constants (Theorem 1). Theorem 1 and a new construction in the proof of Theorem 1 are the  main results of the present paper. Further, we introduce  a generalized  $ K-$functional  which is related to the new approach to the  direct theorems of approximation theory \cite{fks09,bdk12} and  give the analogue of Theorem 1 for it (Theorem 2). We also give a proof of the estimates of the type \eqref{aa} which hold for $\alpha>0$ and better than \eqref{aa} for $\alpha<0.778$ (Theorem 3).

\section{\bf Some auxiliary results }
\vskip0.2cm
\textit{2.1. The classical moduli of continuity and the special moduli of continuity.}
\vskip0.2cm

In this paper, we consider the modulus of continuity of the second order (modulus of smoothness). The classic definition of the modulus of smoothness is the following \cite{dl93}:
\begin{equation*}
\omega_2 (f, h):= \sup_{0<u\leq h}  \sup_{t \in \mathbb{T}} | f(t+u)- 2 f(t) + f(t-u)| =
\sup_{0<u\leq h} \| \Delta_u^2 f \|.
\end{equation*}
In \cite{fks09,bdk12,bks13} the importance of the following moduli of continuity
was indicated
\begin{align*}
&W_2 (f, \chi_h^k)=\|f - f*\chi_h^k \|,  \\
&W_2^* (f, \chi_h^k)=\sup_{0<u\leq h}\|f - f*\chi_u^k \|.
\end{align*}
It is obvious that $W_2 (f, \chi_h^k)\leq W_2^* (f, \chi_h^k)$. Note that
\begin{equation*}
f(t)-(f\ast\chi_h^k)(t)=f(t)-\int_{kh/2}^{kh/2}f(t-u)\chi_h^k(u)\,du
\end{equation*}
\begin{equation*}
=f(t)-\int_0^{kh/2}(f(t+u)+f(t-u))\chi_h^k(u)\,du
=-\int_0^{kh/2}\Delta_u^2f(t)\chi_h^k(u)\,du,
\end{equation*}
and hence
\begin{equation*}
\|f - f*\chi_h^k \| = \| \int_0^{kh/2} \Delta_u^2 f (\cdot) \, \chi_h^k(u)du\|.
\end{equation*}

We give some simple properties of the modulus $ W_2 $.

{\bf Lemma 1.}\ {\it 
Let $f, D^2g \in C(\mathbb{T}),$ $\ k,l \in \mathbb{N}, \ h>0$. Then for $W_2,$ the following inequalities hold
\begin{align}
&W_2 (f, \chi_h^k) \le W_2 (f-g, \chi_h^k) + W_2(g, \chi_h^k), \label{w1}\\
&W_2 (f, \chi_h^k) \le 1/2 \,\omega_2 (f, kh/2),\label{w2}\\
&W_2 (f, \chi_h^k) \le  2 \| f \|, \label{w3} \\
&W_2(f,\chi_h^k)\le W_2(f,\chi_h^l)+W_2(f,\chi_h^{|l-k|}),\quad \mbox{in\, particular}\label{w4}\\
&W_2 (f, \chi_h^{kl}) \le k \, W_2 (f, \chi_h^l),\label{w6}\\
&W_2 (g, \chi_h^k) \le c_k(h)\| D^2 g \|,\ 
 c_k(h)= \frac{2h^2}{(k+2)!}\sum_{j=0}^{[k/2]}(-1)^j \binom kj\left(\frac {k}{2} - j\right)^{k+2}, \label{w5}\\
&\mbox{Specifically, for $k = 1, 2$ \eqref{w5} gives } \nonumber\\ 
&W_2(g, \chi_h) \le \frac{h^2}{24} \| D^2 g  \|, \label{wd}\\
&W_2 (g, \chi_h^2) \le \frac{h^2}{12} \| D^2 g \|. \label{wd2}
 \end{align}
 }

\textbf{Proof.}
The inequalities \eqref{w1},\eqref{w2} follow directly from
the definition of a special modulus of continuity. We have
\begin{equation*}
W_2(f, \chi_h^k) 
= \| (f-g) - (f-g)*\chi_h^k + g - g*\chi_h^k \|
\le  W_2 (f-g, \chi_h^k) + W_2 (g, \chi_h^k)
\end{equation*}
and
\begin{equation*}
W_2 (f, \chi_h^k) = \| \int_0^{kh/2} \Delta_u^2 f(\cdot) \, \chi_h^k (u) \, du \| \le 1/2 \, \omega_2 (f, kh/2).
\end{equation*}
To prove \eqref{w3} it is enough to use the semi-additive property of the norm, the following property of convolution:
$ \| f* \chi_h^k \| \le \|f\|\ast\| \chi_h^k \|_1$
and the equality
$ \|\chi_h^k\|_1 = 1. $
To prove \eqref{w4},  we write  $f - f*\chi_h^{k} $ in the form
\begin{equation*}
f- f*\chi_h^{k} = (f- f*\chi_h^l)+(f*\chi_h^l-f*\chi_h^k),
\end{equation*}
and  similarly to the proof of \eqref{w3}, obtain
\begin{equation*}
W_2(f, \chi_h^{k}) \le  W_2(f,\chi_h^l)+W_2(f,\chi_h^{|l-k|}).
\end{equation*}
To prove \eqref{w5} we will use the representation (see \cite{k}, p. 245)
\begin{equation*}
\chi_h^k (u-kh/2)= \frac 1{h(k-1)!}\sum_{j: j\ge 0, u/h-j>0}  (-1)^j \binom kj \, (u/h-j)^{k-1}, \qquad 0<u<kh.
\end{equation*}
Using inequality $ \|u^{-2}\Delta_u^2 g(\cdot)\|\le \| D^2 g \| $, we obtain
\begin{equation*}
W_2(g, \chi_h^k) \le \| D^2 g \| \int_0^{kh/2} u^2 \chi_h^k (u) \, du 
=\| D^2 g \| \int_0^{kh/2} (kh/2 - u)^2 \chi_h^k (kh/2 - u) \, du 
\end{equation*}
\begin{equation*}
=
\frac{\|D^2 g \|}{h^k (k-1)!} \sum_{j=0}^{[k/2]} (-1)^j \binom kj \int_{jh}^{kh/2}
(u-kh/2)^2 (u-jh)^{k-1}\, du 
\end{equation*}
\begin{equation*}
=\left(\frac{2 h^2}{(k+2)!} \sum_{j=0}^{[k/2]} (-1)^j \binom kj (k/2 -j)^{k+2}\right) \, \| D^2 g \|.\hskip 1cm \square
\end{equation*}


{\bf Remark 1.}\ {\it 
Lemma 1 holds for $ W_2^*$ with the same proofs.
}
\vskip0.2cm
\textit{2.2. The second modulus of continuity and the $K-$functional }
\vskip0.5cm
For $ f \in C (\mathbb{T}) $ define the second $ K-$functional as follows:
\begin{equation*}
K_2 (f, h):=\inf_{g \in C^2 (\mathbb{T})} \{ \|f - g \| + h^2 \| D^2 g \| \}.
\end{equation*}

The second $K-$functional characterizes the  values of the best approximation $f\in C(\mathbb{T})$ by smooth functions $g\in C^2 (\mathbb{T)}, $ with  a  control  on the norm $g$. Note that $K_2$ has the following properties:
\begin{equation*}
K_2 (f,h) \le \max \{1, h^2/\delta^2\} K_2 (f,\delta), \quad h,\delta >0.
\end{equation*}
Indeed, for $\delta\geq h $, we have $\max\{1,h^2/\delta^2\}=1$ and
\begin{equation}\label{k1}
K_2(f,h) \le K_2 (f, \delta).
\end{equation}
If $\delta<h$, then $\max\{1,h^2/\delta^2\}=h^2/\delta^2$ and
\begin{equation*}
h^2/\delta^2 \ K_2(f, \delta)=\inf_{g \in C^2} \{ h^2/\delta^2 \| f - g \| + h^2 \| D^2 g \|\}\ge K_2 (f,h).
\end{equation*}
Thus, when $\delta<h$, we have
\begin{equation}\label{k2}
\frac{K_2(f,h)}{h^2} \le \frac{K_2(f,\delta)}{\delta^2}.
\end{equation}

In other words, the function $K_2(f,h)$ is a monotonically increasing function of the argument $h>0$, and the function $K_2(f,h)/h^2$ is a monotonically decreasing function for $h>0$.

The following lemma is well known. The idea of using an intermediate approximation of a smooth function belongs to Steklov \cite{st24,st83} and Favard \cite{f37}.
To present a standard approach and compare it with our approach, we give the lemma   with the proof.  
Of special interest is the constant in the first inequality in \eqref{l1}.
\vskip0.2cm
{\bf Lemma 2.}\ {\it 
For $f\in C(\mathbb{T}) $ and $h>0$, we have
\begin{equation}\label{l1}
\frac 23 K_2 \, (f, h/2) \le \frac 12 \, \omega_2 (f, h) \le 2 \, K_2 (f, h/2).
\end{equation}
}

\textbf{Proof.}
To prove the right inequality, we will use the well-known properties of the second modulus of continuity
\begin{align*}
\omega_2(f_1+f_2,h) \le \omega_2(f_1,h) + \omega_2(f_2,h),\quad 
\omega_2 (f,h) \le 4 \| f \|, \quad
\omega_2(g,h) \le h^2 \|D^2 g \|.
\end{align*}
We have
\begin{align*}
\omega_2(f,h) &\le \inf_{g \in C^2} \{ \omega_2 (f-g,h) + \omega_2 (g,h) \}\\
&\le \inf_{g \in C^2} \{ 4 \| f - g \| + h^2 \| D^2 g \| \} = 4 \, K_2 (f,h/2).
\end{align*}
If  $ g = f*\chi_h^2$, then by using the identity 
$ D^2(f*\chi_h^2) = h^{-2} \Delta_h^2 f$ we obtain the left inequality
\begin{align*}
K_2(f,h/2) &\le \| f - f*\chi_h^2\| + \frac{h^2}4 \, \| D^2 (f*\chi_h^2) \, \| \\
&\le \frac 12 \, \omega_2 (f,h) + \frac 14 \, \omega_2(f,h) = \frac 34 \, \omega_2(f,h).
\hskip 1cm \square
\end{align*}
\vskip0.2cm
\section{\bf The special modulus of continuity and the $ K-$functional.}

 \vskip0.2cm
 In this and the following sections, we restrict ourselves to statements about the special
moduli of continuity $ W_2 (f,\chi_h^k)$, $ W_2^* (f,\chi_h^k) $ for $k = 1,2. $

The main result of this paper is the following theorem and the new construction \eqref{tau}, which was used in the proof.

\vskip0.2cm
{\bf Theorem 1.}\ {\it
Let $f\in C(\mathbb{T})$ and $h>0$. Then
\begin{equation}\label{t12}
\frac23K_2(f,h/(2\sqrt{6}))\leq W_2^*(f,\chi_{h}^2)\leq 2K_2(f,h/(2\sqrt{6})),
\end{equation}
\begin{equation}\label{t11}
\frac25K_2(f,h/(4\sqrt{3}))\leq W_2^*(f,\chi_{h})\leq 2K_2(f,h/(4\sqrt{3})).
\end{equation}
}

\textbf{Proof.}
Firstly, we prove \eqref{t12}. Let $f\in C(\mathbb{T})$, $g  \in C^2 (\mathbb T)$, $h>0$. The proof of the right inequality is standard. 
By \eqref{w1}, \eqref{w3}, \eqref{wd2} we have:
\begin{align*}
W^*_2(f,\chi_h^2)&\leq\inf_{g\in C^2}\{W^*_2(f-g,\chi_h^2)+W^*_2(g,\chi_h^2)\}\\
&\leq \inf_{g\in C^2}\{2\|f-g\|+\frac{h^2}{12}\|D^2 g\|\}=2 K_2(f,h/(2\sqrt{6})).
\end{align*}
Consider the inverse estimate. Let
\begin{equation}\label{tau}
g(t)=\frac{12}{h^2}\int_0^{h}(f\ast\chi_u^2)(t)\,u^2\chi_{h}^2(u)\, du.
\end{equation}
From the definition of the second $K-$functional we obtain
\begin{equation*}
K_2(f,h/(2\sqrt{6}))
\leq\|f-g\|+\frac{h^2}{24}\|D^{2}g\|.
\end{equation*}
We can compute 
\begin{equation}\label{D^2 tau}
\|D^2g\|=\frac{12}{h^2}W_2(f,\chi_{h}^2)\leq\frac{12}{h^2}W^*_2(f,\chi_{h}^2).
\end{equation}
Indeed,  since  $D^2(f\ast\chi_u^2)={u^{-2}}\Delta_u^2f,$ we have
\begin{align*}
D^2g(t)
&=\frac{12}{h^2}\int_0^{h}D^2((f\ast\chi_u^2)(t))\,u^2\chi_{h}^2(u)\, du\\
&=\frac{12}{h^2}\int_0^{h}\Delta_u^2f(t)\,\chi_{h}^2(u)\, du
=-\frac{12}{h^2}\left(f(t)-(f\ast\chi_{h}^2)(t)\right),
\end{align*}
which implies \eqref{D^2 tau}.
Further, from ${12}/{h^2}\int_0^{h} u^2\chi_{h}^2(u)\, du =1$ (see the proof of \eqref{w5} and \eqref{wd2}), we have
\begin{equation}\label{f-ta}
\begin{aligned}
|f(t)-g(t)|
&=\left|\frac{12}{h^2}\int_0^{h}\left(f(t)-(f\ast\chi_u^2)(t)\right)u^2\chi_{h}^2(u)\, du\,\right| \\
&\leq W^*_2(f,\chi_{h}^2)\ \frac{12}{h^2}\,\int_0^{h}u^2\chi_{h}^2(u)\, du
= W^*_2(f,\chi_{h}^2). 
\end{aligned}
\end{equation}
From \eqref{D^2 tau} and \eqref{f-ta} we obtain
\begin{equation*}
K_2(f,h/(2\sqrt{6}))
\leq W^*_2(f,\chi_{h}^2)+\frac{1}{2}W^*_2(f,\chi_{h}^2).
\end{equation*}
We proceed with the proof of \eqref{t11}. The right inequality 
of \eqref{t11} can be proved as in the previous case.
\begin{align*}
W_2^*(f,\chi_{h})&\leq \inf_{g\in C^2}\{W_2^*(f-g,\chi_{h})+W_2^*(g,\chi_{h})\}\\
&\leq \inf_{g\in C^2}\{2\|f-g\|+\frac{h^2}{24}\|D^2 g\|\}=2K_2(f,h/(4\sqrt{3})).
\end{align*}

In the proof of the inverse estimate we will use two auxiliary functions, 
$g_1(t)$ and $g_2(t)$.
\begin{equation*}
g_1(t)=\frac{24}{h^2}\int_0^{h/2}(f\ast\chi_u)(t)\,u^2\chi_{h}(u)\, du,
\end{equation*}
\begin{equation*}
g_2(x)=\frac{24}{h^2}\int_0^{h/2}(f\ast\chi_u^2)(t)\,u^2\chi_{h}(u)\, du.
\end{equation*}
From the definition of the second $K-$functional it follows that
\begin{equation}\label{K_2}
K_2(f,h/(4\sqrt{3}))
\leq \|f-g_1\|+\|g_1-g_2\|+\frac{h^2}{48}\|D^{2}g_2\|.
\end{equation}
Similarly, \eqref{D^2 tau} and \eqref{f-ta} yield
\begin{align}
\|f-g_1\|&\leq W^*_2(f,\chi_{h}), \label{f-g_1}\\
\|D^2g_2\|&\leq\frac{24}{h^2}\,W^*_2(f,\chi_{h}).\label{D^2 g_2}
\end{align}
It remains to estimate the norm $ g_1 - g_2$. We have
\begin{equation*}
|g_1(t)-g_2(t)|
=\left|\frac{24}{h^2}\int_0^{h/2}\left((f\ast\chi_u)(t)-(f\ast\chi_u^2)(t)\right)u^2\chi_{h}(u)\, du \, \right|.
\end{equation*}
By using the inequalities $\|f\ast g\|\leq\|f \|\|g \|_1$ and $ W^*_2(f,\chi_{h/2})\leq W^*_2(f,\chi_{h})$ and \eqref{wd}, we obtain
\begin{equation*}
\|g_1-g_2\|
\leq W^*_2(f,\chi_{h})\ \frac{24}{h^2} \, \int_0^{h/2}u^2\chi_{h}(u)\, du
= W^*_2(f,\chi_{h}).
\end{equation*}
By \eqref{K_2},  \eqref{f-g_1}, \eqref{D^2 g_2} and the previous inequality, we get
\begin{equation*}
K_2(f,h/(4\sqrt{3}))
\leq
2W^*_2(f,\chi_{h})+\frac{1}{2}W^*_2(f,\chi_{h}).\hskip 1cm \square
\end{equation*}

Theorem $1$ and Lemma $2$ yield the following fact.
\vskip .2cm
{\bf Corollary 1.}\ {\it
 Let $h>0$. Then for $h_1=h/(2\sqrt{6})$, $h_2=h/(4\sqrt{3})$
\begin{equation*}
\frac13 \, W_2^*(f,\chi_h^2)\leq\frac23K_2(f,h_1)
\leq\frac12\omega_2(f,2h_1)\leq 2K_2(f,h_1)\leq3W_2^*(f,\chi_h^2).
\end{equation*}
and
\begin{equation*}
\frac13 W_2^*(f,\chi_h)\leq\frac23K_2(f,h_2)
\leq\frac12\omega_2(f,2h_2)
\leq2K_2(f,h_2)
\leq 5W_2^*(f,\chi_h).
\end{equation*}
}

Moreover, Theorem 1 and the inequality \eqref{w4} of Lemma 1 allow us to compare the moduli $W_2^*(f,\chi_h)$ and $W_2^*(f,\chi_h^2)$:

{\bf Corollary 2.} {\it
For  $h>0$
\begin{equation}\label{sl2}
W_2^*(f,\chi_h^2) \le 2 \, W_2^* (f, \chi_h) \le 6 \, W_2^*(f,\chi_h^2),
\end{equation}
where the left inequality cannot be improved for $h= 1/(2n), \ \  n=1,2, \dots.$
}

\textbf{Proof.}
The left inequality in \eqref{sl2} is \eqref{w4} for $k=2,\ l=1$. The right inequality follows from  \eqref{t11} and  \eqref{k1}:
\begin{equation*}
W_2^*(f,\chi_h) \le 2 \, K_2 (f, h/(4\sqrt{3})) \le 2 \, K_2 (f, h/(2\sqrt{6})) \le 3 \, W_2^* (f,\chi_h^2).
\end{equation*}

The left inequality of Corollary 2 is exact for $h=1/(2n),\ n=1,2, \dots$.
To show this, fix $ n \in \mathbb{N} $ and $ h=1/(2n)$. Consider the following construction, which was used in \cite {bdk12}. Let
\begin{align*}
\varepsilon_j(t)&:= \mbox{sign\,} \cos (2\pi n t)* \chi_h^j(t), \qquad j=0,1,2,\dots, \quad \chi_h^0= \delta, \quad f*\delta = f, \\
\phi &:= \sum_{jем=0}^\infty \varepsilon_j.
\end{align*}

It was proved in \cite{bks13} that $\|\varepsilon_j\| = E_{n-1} (\chi_h^j) \le (2 / \pi)^{j-1}$ and  the series converges uniformly. We have
\begin{align*}
\| \phi - \phi*\chi_h^2 \| &= \| \sum_{j=0}^\infty \varepsilon_j - \sum_{j=2}^\infty \varepsilon_j \| = \| \varepsilon_0 + \varepsilon_1 \| =2, \\
\|\phi - \phi*\chi_h\| &= \| \sum_{j=0}^\infty \varepsilon_j - \sum_{j=1}^\infty \varepsilon_j \|=
\| \varepsilon_0 \| =1.\hskip 1cm \square
\end{align*}

Now we show that Theorem 1 implies  the estimates  \eqref{k}, obtained in  \cite{bdk12} by another method.

{\bf Corollary 3.} {\it
For  $h>0$
\begin{equation*}
\frac 25 K_2 (f, h/(4\sqrt{6})) \le 
W_2^*(f,\chi_h) \le 4 \, K_2 (f, h/(4\sqrt{6})).
\end{equation*}
}

\textbf{Proof.}
The fact that $K_2(h)$ is increasing  plus the left inequality \eqref{t11} imply
\begin{equation*}
\frac 25 \, K_2 (f, h/(4\sqrt{6})) \le \frac 25 \, K_2(f,h/(4\sqrt{3})) \le W_2^*(f,\chi_h).
\end{equation*}
The property \eqref{k2} and the right inequality in \eqref{t11} give
\begin{equation*}
W_2^*(f,\chi_h) \le 2 \, K_2 (f, h/(4\sqrt{3})) \le 4 \, K_2 (f, h/(4\sqrt{6})).\hskip 1cm \square
\end{equation*}


\vskip0.2cm
Note that the constant in the lower estimate \eqref{t11} is worse than the constant in the appropriate lower estimate \eqref{t12}. The estimate \eqref{t12} can be improved if we consider the generalized $\widetilde K_2$-functional, which is a sharper characteristic than
$K_2$. Its definition is motivated by the proof of \, \eqref{t11}. 
\vskip .2cm
Let 
\begin{equation*}
\widetilde{K}_2(f,h_1,h_2)=\inf_{g_1\in C,g_2\in C^2 }\{\|f-g_1\|+h_1\|D(g_1-g_2)\|+h_2^2\|D^2g_2\|\}.
\end{equation*}
Clearly, one can take $ g_2\equiv g_1$ and obtain
\begin{equation*}
\widetilde{K}_2(f,h_1,h_2) \le K_2(f,h_2).
\end{equation*}
{\bf Theorem 2.} {\it
 Let $f\in C(\mathbb{T})$ and $h>0$. Then
\begin{equation*}
\frac{4}{9}\widetilde{K}_2\left(f,\frac{h}{8},\frac{h}{4\sqrt{3}}\right)\leq W_2^{*}(f,h)\leq2\widetilde{K}_2\left(f,\frac{h}{8},\frac{h}{4\sqrt{3}}\right).
\end{equation*}
}

\textbf{Proof.}
Let $f, D g_1, D^2g_2 \in C(\mathbb{T})$ and $h>0.$ 
By \eqref{w1}, \eqref{w3}, \eqref{wd} and the inequality
\begin{align*}
W_2^{*}(g_1-g_2,h)&\leq\|\int_{0}^{h/2}\Delta_u^1\Delta_u^1(g_1-g_2)(\cdot)\chi_h(u)\,du \| \\
&\leq \|D(g_1-g_2)\|\, \frac{2}{h} \,\int_{0}^{h/2}u\,du
=\frac{h}{4}\|D(g_1-g_2)\|
\end{align*}
\noindent
we obtain
\begin{equation*}
W_2^{*}(f,h)\leq W_2^{*}(f-g_1,h)+W_2^{*}(g_1-g_2,h)+W_2^{*}(g_2,h)
\leq 2\widetilde{K}_2\left(f,\frac{h}{8},\frac{h}{4\sqrt{3}}\right).
\end{equation*}
Consider the inverse estimate. By the definition of $\widetilde{K}_2(f,h_1,h_2) $, we have
\begin{equation}\label{sK2}
\widetilde{K}_2\left(f,\frac{h}{8},\frac{h}{4\sqrt{3}}\right)
\leq\|f-g_1\|+\frac{h}{8}\|D(g_1-g_2)\|+\frac{h^2}{48}\|D^2 g_2 \|,
\end{equation}
where
\begin{equation*}\label{varphi_2}
g_1(t)=\frac{24}{h^2}\int_0^{h/2}(f\ast\chi_u)(t)\,u^2\chi_{h}(u)\, du,
\end{equation*}
\begin{equation*}
g_2(t)=\frac{24}{h^2}\int_0^{h/2}(f\ast\chi_u^2)(t)\,u^2\chi_{h}(u)\, du.
\end{equation*}
We can estimate the norm of $D(g_1-g_2).$  From the inequality
\begin{equation*}
\|D(f\ast\chi_u-f\ast\chi_u^2)\|
\leq 2u^{-1}W_2^*(f,\chi_u)
\end{equation*}
we conclude that
\begin{equation*}
\|D(g_1-g_2)\|
\leq
\frac{48}{h^2}\, W_2^*(f,\chi_h)\int_0^{h/2} \, u\,\chi_{h}(u)\, du=
\frac{6}{h}W_2^*(f,\chi_h).
\end{equation*}
Therefore, by \eqref{sK2}, \eqref{f-g_1}, \eqref{D^2 g_2}
and the previous inequality, we deduce
\begin{equation*}
\widetilde{K}_2\left(f,\frac{h}{8},\frac{h}{4\sqrt{3}}\right)
\leq W_2^*(f,\chi_h)+\frac{3}{4} W_2^*(f,\chi_h)+\frac{1}{2}W_2^*(f,\chi_h)
=\frac{9}{4}W_2^*(f,\chi_h).\hskip 1cm \square
\end{equation*}

\vskip0.2cm
\section{\bf Jackson's theorem with the special modulus of continuity} 

 To prove Jackson's theorem, we will use a method that dates back to  Steklov\cite{st83, st24}  and Favard \cite{f37}.
 A new element here is our intermediate function $g$.
 
{\bf Theorem 3.} {\it
Let $f\in C(\mathbb{T})$, $\alpha>0$, $c(\alpha)=1+{3}/(2\alpha^2)$. Then for $n\in N$
\begin{equation}\label{j1}
E_{n-1}(f) \leq c(\alpha)\, W_2^{*}(f,\chi_{\alpha/(2n)}^2),
\end{equation}
\begin{equation}\label{j2}
E_{n-1}(f) \leq 2c(\alpha)\, W_2^{*}(f,\chi_{\alpha/(2n)}).
\end{equation}
}

\textbf{Proof.}
Firstly, we prove the inequality \eqref{j1}. 
As an intermediate approximation we will use

\begin{equation*}
g(t)=\frac{12}{h^2}\int_0^{h}(f\ast\chi_u^2)(t)\,t^2\chi_{h}^2(u)\, du,\quad h>0.
\end{equation*}
We can write the function $ f $ as

\begin{equation*}
f=f-g+g,
\end{equation*}
and obtain
\begin{equation*}
E_{n-1}(f) \leq\|f-g\|+E_{n-1}(g).
\end{equation*} 
We have (see \eqref{f-ta})

\begin{equation*}
\|f-g\|\leq W_2^{*}(f,\chi_{h}^2).
\end{equation*}
By using the Favard inequality  
\begin{equation}\label{akf}
E_{n-1}(g) \le \frac{1}{32n^2} \|D^2 g \|
\end{equation}
and \eqref{D^2 tau}, we see that
\begin{equation*}
E_{n-1}(g) \leq\frac{1}{32n^2}\|D^2g\|
\leq
\frac{3}{8(nh)^2}W_2^{*}(f,\chi_{h}^2).
\end{equation*}
By choosing $ h = \alpha/(2n) $, $\alpha>0$, $n\in \mathbb{N}$ we finally obtain 
\begin{equation*}
E_{n-1}(f) \leq\left(1+\frac{3}{2\alpha^2}\right)W_2^{*}(f,\chi_{\alpha/(2n)}^2).
\end{equation*}
The inequality \eqref{j2} follows from Corollary 2 and the inequality \eqref{j1}. \hskip 1cm $\square$
\vskip0.2cm
The estimates of Theorem 3 are valid for $\alpha >0$. The inequality \eqref{j2}  complements the inequality
\eqref{aa}, which is valid for $\alpha > 2/\pi$. Note that the constant in \eqref{j2}  is better than in \eqref{aa} for $ \alpha \leq 0.778.$ 

	For the classical modulus of smoothness $\omega_2$ the standard Steklov--Favard approach
(see for example \cite[Theorem 9.2]{fks09})
gives  
\begin{equation}\label{zs}
E_{n-1}(f) \le \left( \frac 12 + \frac{1}{8\alpha^2} \right) \, 
\omega_2 \left(f, \frac{\alpha}{2n} \right), \qquad \alpha >0.
\end{equation}

The inequality \eqref{zs} is sharp for $\alpha = 1/(2j), \ j=1,2, \dots $ (see \cite{p13}).
	The estimates \eqref{zs} and $ \omega_2 (f, h) \le h^2 \| D^2 f \|$ imply the
weak form of the Favard inequality,
$$
E_{n-1}(f) \le \frac 1{32n^2} \left( 1 + 4 \alpha^2 \right) \, \|D^2 f \|.
$$
The value $4 \alpha^2$ is a measure of the distance between the sharp difference estimate
\eqref{zs} and the sharp differential
estimate \eqref{akf}.

By the inequalities
$$
W_2^*(f,\chi_h)\leq\frac{h^2}{24}\|D^2f\|,\quad W_2^*(f,\chi_h^2)\leq\frac{h^2}{12}\|D^2f\|,
$$
and  \eqref{j1}, \eqref{j2} we have

$$
E_{n-1}(f) \leq \frac{1}{32n^2}\left(1 + \frac{2\alpha^2}{3} \right)\|D^2f\|.
$$

 Probably, by using sharp estimates in Theorem 3, for small $ \alpha, $
one can reduce the $\alpha^2$-coefficient 2/3  and obtain the  better correspondence between
differential  and difference estimates  in the direct theorems of approximation theory.

\section*{Acknowledgment}

The author is greatly indebted to Y. Kryakin for suggesting the problem and for many
stimulating conversations and valuable comments.

\bibliographystyle{model1a-num-names}
\bibliography{<your-bib-database>}



\end{document}